\documentclass{article}
\usepackage{epsf}
\usepackage{fullpage}
\usepackage{euscript,amsmath,amssymb, amsthm} 
\usepackage[dvips]{graphicx}
\usepackage[dvips]{color}
\renewcommand{\k}{\mathbf{k}}
\newcommand{\C}{\mathbf{C}}
\newcommand{\Z}{\mathbf{Z}}

\newcommand{\U}{\mathfrak{U}}

\renewcommand{\span}{\operatorname{\mathsf{span}}}
\renewcommand{\mod}{\mathsf{\ mod \,}}
\renewcommand{\dim}{\operatorname{\mathsf{dim}}}

\newcommand{\osp}{\operatorname{\mathfrak{osp}}}

\renewcommand{\a}{\mathsf{a}}

\newcommand{\F}{\mathsf{F}}

\newcommand{\s}{\mathsf{s}}
\newcommand{\x}{\mathsf{x}}
\newcommand{\y}{\mathsf{y}}

\newcommand{\isom}{\cong}

\newcommand{\impl}{\Longrightarrow}

\newcommand{\defeq}{\stackrel{\mathrm{def}}{=}}

\def\cent#1{\mathcal{Z}(#1)}

\newcommand{\Hh}{\mathbf{H}}
\newcommand{\hhh}{\mathsf{h}}
\newcommand{\ert}{\varepsilon}
\newcommand{\pf}{\mathbf{F}_p}

\begin{document}

\newtheorem{theorem}{Theorem}[section]
\newtheorem{proposition}[theorem]{Proposition}
\newtheorem{lemma}[theorem]{Lemma}
\newtheorem{claim}[theorem]{Claim}  
\newtheorem{corollary}[theorem]{Corollary}

\theoremstyle{definition}
\newtheorem{definition}[theorem]{Definition}
\newtheorem{definition-sort-of}[theorem]{`Definition'}
\newtheorem{notation}[theorem]{Notation}

\newtheorem{example}[theorem]{Example}
\newtheorem{remark}[theorem]{Remark}
\newtheorem{fact}[theorem]{Fact} 
\newenvironment{block}[1]{\medskip \noindent {\Large \bf #1.}}{\medskip}

\title{Representations of rational Cherednik algebras of rank 1 
in positive characteristic}

\author{Fr\'ed\'eric Latour}

\maketitle

\section{Introduction}

Let $\Gamma\subset SL_2(\C)$ be a finite subgroup, 
and $r$ be the number of conjugacy classes of $\Gamma$.
Let $D$ be the algebra of polynomial differential operators 
in one variable with complex coefficients. 
In the paper \cite{CBH}, Crawley-Boevey and Holland 
introduced an $r$-parameter 
family of algebras $\mathbf H_{t,c_1,...,c_{r-1}}$ over $\C$, 
which is a universal deformation of 
the semidirect product algebra $\C\Gamma\ltimes D$. They also 
described their representation theory, showing that it is related 
to the structure of the root system attached to $\Gamma$ 
via the McKay's correspondence. 

The algebras of Crawley-Boevey and Holland are the simplest case 
of so called symplectic reflection algebras introduced 
in \cite{EG}. These algebras have interesting representation theory,
about which little is known. The better understood case is that of 
rational Cherednik algebras corresponding to complex reflection groups
\cite{DO,GGOR,CE}, and especially that of real reflection 
(in particular, Weyl) groups \cite{BEG,G}. In the Weyl
case, these algebras admit a q-deformation
defined by Cherednik, 
whose representation theory (while having some similarities
to the rational case) 
can be effectively studied by the technique of intertwiners \cite{Ch}
or geometric representation theory \cite{Va}. 

The goal of this paper is to begin the study of representation theory 
of rational Cherednik algebras in positive characteristic. 
Namely, we study the simplest case of rank 1,
i.e. the algebras of 
\cite{CBH} for $\Gamma=\mathbf Z/r\mathbf Z$ in characteristic $p$, 
which is prime to $r$. Our main result 
is a complete description of irreducible representations 
of such algebras. 

The paper is organized as follows. 

In Section 2, we state the main results. 

In Section 3, we prove the results for the ``quantum'' case, i.e. 
the case when the ``Planck's constant'' $t$ is nonzero. In this case, 
generic irreducible representations have dimension
$pr$. 

In Section 4, 
we prove the results for the ``classical'' case, i.e. 
the case when the ``Planck's constant'' $t$ is zero. In this case, 
generic irreducible representations have dimension
$r$. 

{\bf Acknowledgements.} The author thanks his adviser Pavel Etingof for
posing the problem and useful discussions, as well as for helping to write 
an introduction. The work of the author 
was partially supported by the National Science Foundation (NSF) grant 
DMS-9988796 and by a Natural Sciences and Engineering Research 
Council of Canada (NSERC) Julie Payette research scholarship. 

\section{Statement of results}

Let $\k$ be an algebraically closed field of characteristic $p$ and $r \geq 1$ 
a positive integer relatively prime to $p$. Let $\ert \in \k$ be a 
primitive $r$-th root of unity.

\begin{definition} Let $t,c_1,...,c_{r-1}\in \k$. 
The algebra $\Hh_{t,c_1,...,c_{r-1}}$ over $\k$ 
is generated by $\x,\y,\s$ with relations
$$\s^r = 1, \quad \s\x = \ert^{-1} \x\s, \quad \s\y = \ert\y\s, \quad 
[\y,\x] = t - \sum_{j=1}^{r-1} c_j \s^j.$$
\end{definition}

In characteristic zero, this algebra has been studied in \cite{CBH}.
In the sequel, we fix $c_1,...,c_{r-1}$ and denote the algebra 
$\Hh_{t,c_1,...,c_{r-1}}$ simply by $\Hh_t$. 

We will classify irreducible representations of $\Hh_t$. 
In particular, we show that they are all finite dimensional. 

It is clear that for $t\ne 0$, $\Hh_t$ is isomorphic to $\Hh_1$ 
with rescaled parameters $c_j$. Thus it is sufficient to 
classify irreducible representations of $\Hh_0$ and $\Hh_1$. So we set 
$\Hh_1=\Hh$, $\Hh_0=\Hh'$. 

\begin{remark}
In the trivial case $r=1$, $\Hh'$ is the polynomial algebra in 
two variables and $\Hh$ is the Weyl algebra. In the case $r=2$,
as pointed out in \cite{MC},
$\Hh$ is a quotient of the universal enveloping algebra 
$\U(\osp(1,2))$ of the Lie superalgebra $\osp(1,2)$
by a central character.
\end{remark}

Below we will need the element
$$
\hhh \defeq \x\y -\sum_{j=1}^{r-1} c_j (1 - \ert^{-j})^{-1} \s^j,
$$
both in $\Hh$ and $\Hh'$. This element (in the case $r=2$) was 
introduced in \cite{MC}. 

\subsection{Irreducible representations of $\Hh$} 

The elements $\x^{pr},\y^{pr}$ are central in $\Hh$. Thus they have to act as
 scalars in irreducible representations. 

\begin{proposition}
Suppose that $a,b \in \k$ are not both zero.
Then, for each $\beta\in \k$ which satisfies the equation
\begin{equation}
\prod_{m=0}^{r-1}
f \left(\beta+\sum c_j(1-\varepsilon^{-j})^{-1}\varepsilon^{mj}\right)=ab,
\label{propeqn1} 
\end{equation}
where $f(z)=z^p-z$, there exists a unique irreducible
representation $V_{\beta,a,b}$ of $\Hh$, where
$\x^{pr}$ acts on $V_{\beta,a,b}$ through multiplication by $a$, $\y^{pr}$
acts on $V_{\beta,a,b}$ through multiplication by $b$, 
and the spectrum of $\hhh$ 
in $V_{\beta,a,b}$ is $\beta+\mathbf F_p$.

Explicitly, the representation
$V_{\beta,a,b}$ looks as follows. 
Let $$\mu_i = \beta + i + \sum_{j=1}^{r-1} c_j (1 - \ert^{-j})^{-1} \ert^{-ij},$$ 
for all $i=0,...,pr-1$.
Then for $a \neq 0$ the representation $V_{\beta,a,b}$
 has a basis $\{ v_0, v_1, \ldots, v_{pr-1} \}$
in which
\begin{eqnarray*}
\x v_i &=& v_{i+1} \quad \mbox{for } i = 0, \ldots, pr-2, \\
\x v_{pr-1} &=& a v_0, \\
\y v_i &=& \mu_i v_{i-1} \quad \mbox{for } i = 1, \ldots, pr-1, \\
\y v_0 &=& {\it a}^{-1} \mu_0 v_{pr-1}, \\
\s v_i &=& \ert^{-i} v_i. \quad \mbox{for } i = 0, \ldots, pr-1.
\end{eqnarray*}
On the other hand, if $a = 0$ (but $b \neq 0$), we have
\begin{eqnarray*}
\x v_i &=& \mu_{i+1} v_{i+1} \quad \mbox{for } i = 0, \ldots, pr-2, \\
\x v_{pr-1} &=& b^{-1} \mu_0 v_0, \\
\y v_i &=& v_{i-1} \quad \mbox{for } i = 1, \ldots, pr-1, \\
\y v_0 &=& b v_{pr-1}, \\
\s v_i &=& \ert^{-i} v_i. \quad \mbox{for } i = 0, \ldots, pr-1.
\end{eqnarray*}
Thus $\dim V_{\beta,a,b} = pr$.
Two roots $\beta$ and $\beta'$ of \eqref{propeqn1} yield isomorphic   
representations if and only if $\beta - \beta' \in \F_p$.
Furthermore, every irreducible representation of $\Hh$ on which
$\x^{pr}$ and $\y^{pr}$ do not both act as zero is of the form
$V_{\beta,a,b}$
for some $\beta,a,b$.
\label{prop1}
\end{proposition}  

\begin{example} 
For $r=2$, $V_{\beta,a,b} (a \neq 0)$ will has the following 
form:
\begin{equation*} 
\x \mapsto \left(\begin{array}{ccccc} 
0 & 0 & \cdots & 0 & a \\ 
1 & 0 & \cdots & 0 & 0 \\ 
0 & 1 & \cdots & 0 & 0 \\ 
\vdots & \vdots & \ddots & \vdots & \vdots \\
0 & 0 & \cdots & 1 & 0 
\end{array}\right), 
\y \mapsto \left(\begin{array}{ccccc} 
0 & \mu_1 & 0 & \cdots & 0 \\ 
0 & 0 & \mu_2 & \cdots & 0 \\ 
\vdots & \vdots & \vdots & \ddots & \vdots \\ 
0 & 0 & 0 & \cdots & \mu_{2p-1} \\ 
\a^{-1} \mu_0 & 0 & 0 & \cdots & 0 
\end{array}\right), 
\s \mapsto \left(\begin{array}{ccccc} 
1 & 0 & 0 & \cdots & 0 \\ 
0 & -1 & 0 & \cdots & 0 \\ 
0 & 0 & 1 & \cdots & 0 \\
\vdots & \vdots & \vdots & \ddots & \vdots \\
0 & 0 & 0 & \cdots & -1 
\end{array}\right).
\end{equation*} 
\end{example}

\begin{proposition}
For each $m \in \Z/r\Z$, there exists a unique irreducible representation
$W_m$ of $\Hh$ with $\x^{pr}  = \y^{pr}  = 0$ and $\s v = \ert^m v$
whenever $\y v = 0$.
The dimension $D \leq pr$ of $W_m$ is the smallest positive integer
satisfying the equation
\begin{equation*}
D = \sum_{i=0}^{D-1} \sum_{j=1}^{r-1} c_j \ert^{(m-i)j}.
\end{equation*}
$W_m$ has a basis $\{ v_0, v_1, \ldots, v_{D-1} \}$, such that
\begin{eqnarray*}
\x v_i &=& v_{i+1}, \quad \mbox{for } 0 \leq i \leq D-2, \\
\x v_{D-1} &=& 0, \\
\y v_i &=& \mu_i v_{i-1} \quad \mbox{for } 1 \leq i \leq D-1, \\
\y v_0 &=& 0, \\
\s v_i &=& \ert^{m-i} v_i \quad \mbox{for } 0 \leq i \leq D-1,
\end{eqnarray*}  
where $\mu_k = k - \sum_{i=0}^{k-1} \sum_{j=1}^{r-1} c_j \ert^{(m-i)j}.$
Every irreducible representation of $\Hh$ on which
$\x^{pr}$ and $\y^{pr}$ both act as zero is of the form
$W_m$ for some $m$.
\label{prop2}
\end{proposition}   

\begin{remark}
For generic $c_j$, we have $\dim W_m = pr$ for all $m$. 
\end{remark}

\subsection{Irreducible representations of $\Hh'$}

The elements $\x^r,\y^r$ and $\hhh$ are central in $\Hh'$, 
so they must act as scalars in irreducible representations. 

\begin{proposition}
Suppose that $a,b \in \k$ are not both zero.
Then, for each $\beta$ which satisfies the equation
\begin{equation}
ab = \prod_{m=0}^{r-1} \left( \beta + \sum_{j=1}^{r-1} c_j
        (1 - \ert^{-j})^{-1} \ert^{mj} \right),
\label{0propeqn1}
\end{equation}
there exists a unique irreducible
representation $V_{\beta,a,b}$ of $\Hh'$, where 
$\x^{r}$ acts on $V_{\beta,a,b}$ through multiplication by $a$, $\y^{r}$
acts on $V_{\beta,a,b}$ through multiplication by $b$, and 
$\hhh$ acts on $V_{\beta,a,b}$ through 
multiplication by $\beta$.

Explicitly, the representation
$V_{\beta,a,b}$ looks as follows. 
Let $$\mu_i = \beta  + \sum_{j=1}^{r-1} c_j (1 - \ert^{-j})^{-1} \ert^{-ij},$$ 
for all $i=0,...,r-1$.
Then for $a \neq 0$ the representation $V_{\beta,a,b}$
 has a basis $\{ v_0, v_1, \ldots, v_{r-1} \}$
in which
\begin{eqnarray*}
\x v_i &=& v_{i+1} \quad \mbox{for } i = 0, \ldots, r-2, \\
\x v_{r-1} &=& a v_0, \\
\y v_i &=& \mu_i v_{i-1} \quad \mbox{for } i = 1, \ldots, r-1, \\
\y v_0 &=& a^{-1} \mu_0 v_{r-1}, \\
\s v_i &=& \ert^{-i} v_i. \quad \mbox{for } i = 0, \ldots, r-1.
\end{eqnarray*}
On the other hand, if $a = 0$ (but $b \neq 0$), we have
\begin{eqnarray*}  
\x v_i &=& \mu_{i+1} v_{i+1} \quad \mbox{for } i = 0, \ldots, r-2, \\
\x v_{r-1} &=& b^{-1} \mu_0 v_0, \\
\y v_i &=& v_{i-1} \quad \mbox{for } i = 1, \ldots, r-1, \\
\y v_0 &=& b v_{r-1}, \\
\s v_i &=& \ert^{-i} v_i. \quad \mbox{for } i = 0, \ldots, r-1.
\end{eqnarray*} 
Thus $\dim V_{\beta,a,b} = r$.
Every irreducible representation of $\Hh'$ on which
$\x^{r}$ and $\y^{r}$ do not both act as zero is of the form
$V_{\beta,a,b}$ for some $\beta,a,b$.
\label{0prop1}
\end{proposition}

\begin{proposition}
For each $m \in \Z/r\Z$, there exists a unique irreducible representation
$W_m$ of $\Hh'$ with $\x^{r}  = \y^{r}  = 0$ and $\s v = \ert^m v$
whenever $\y v = 0$.
The dimension $D \leq r$ of $W_m$ is the smallest positive integer
satisfying the equation
\begin{equation*}
0 = \sum_{i=0}^{D-1} \sum_{j=1}^{r-1} c_j \ert^{(m-i)j}.
\end{equation*} 
$W_m$ has a basis $\{ v_0, v_1, \ldots, v_{D-1} \}$, such that
\begin{eqnarray*}
\x v_i &=& v_{i+1}, \quad \mbox{for } 0 \leq i \leq D-2, \\
\x v_{D-1} &=& 0, \\
\y v_i &=& \mu_i v_{i-1} \quad \mbox{for } 1 \leq i \leq D-1, \\
\y v_0 &=& 0, \\
\s v_i &=& \ert^{m-i} v_i \quad \mbox{for } 0 \leq i \leq D-1,
\end{eqnarray*}
where $\mu_k = - \sum_{i=0}^{k-1} \sum_{j=1}^{r-1} c_j \ert^{(m-i)j}.$
Every irreducible representation of $\Hh'$ on which
$\x^{r}$ and $\y^{r}$ both act as zero is of the form
$W_m$ for some $m$.
\label{0prop2}
\end{proposition}

\section{Proof of propositions \ref{prop1} and \ref{prop2}}

\begin{lemma}[PBW for $\Hh$, easy direction]
The elements $$\x^i \y^j \s^l, \quad \quad 0 \leq i,j, \quad 0 \leq l 
\leq r-1$$ span $\Hh$ over $\k$.
\label{PBW}
\end{lemma}

\begin{proof}
Given a product of $\x, \y, \s$ in any order, one can ensure that the
$\y$'s are to the right of all the $\x$'s by using $\y\x = \x\y + 1 - 
\sum_{j=1}^{r-1} c_j \s^j$ repeatedly, and one can also ensure that 
the $\s$'s are to the right of all the $\x$'s and $\y$'s by using $\s\x = 
\ert^{-1} \x\s$ and $\s\y = \ert \y\s$ repeatedly. 
\end{proof}

\begin{lemma}
$\x^{pr}$ and $\y^{pr}$ belong to the center $\cent{\Hh}$ of $\Hh$.
\label{lemma1}
\end{lemma}

\begin{proof}
First, let us show that $\x^{pr} \in \cent{\Hh}$.
It is enough to show that $[\x^{pr},\y] = [\x^{r},\s] = 0$. 
But $\x^r \s = \ert^r \s \x^r = \s \x^r$, so $[\x^r, \s] = 0$.
Now let $g(\s) = 1 - \sum_{j=1}^{r-1} c_j \s^j.$ Then,
\begin{eqnarray*}
\y\x^r &=& (\y\x)\x^{r-1} \\
&=& (\x\y + g(\s)) \x^{r-1} \\
&=& \x (\y \x) \x^{r-2} + g(\s) \x^{r-1} \\
&=& \cdots \\
&=& \x^r\y + g(\s) \x^{r-1} + \x g(\s) \x^{r-2} + \cdots + \x^{r-1} g(\s) 
\\
&=& \x^r\y  + r \x^{r-1} + \sum_{j=1}^{r-1} c_j \left( \s^j 
\x^{r-1} + \x \s^j \x^{r-2} + \cdots + \x^{r-1} \s^j \right) \\
&=& \x^r\y + r \x^{r-1} + \sum_{j=1}^{r-1} c_j (1 + \ert^j + 
\cdots + \ert^{j(r-1)}) \s^j \x^{r-1} \\ 
&=& \x^r\y + r \x^{r-1},
\end{eqnarray*}
since $1 + z + \cdots + z^{r-1} = 0$ whenever $z \neq 1$ is an 
$r$-th root of unity. Thus $[\y, \x^r] = r \x^{r-1}$, and hence 
$[\y, \x^{pr}] = pr \x^{pr-1} = 0$.

This completes the proof that $\x^{pr} \in \cent{\Hh}.$ The proof that 
$\y^{pr} \in \cent{\Hh}$ is similar.
\end{proof}

\begin{corollary}
$\Hh$ is finitely generated as a module over its center.
\label{cor2}
\end{corollary}

\begin{proof}
From Lemmas
\ref{PBW} and \ref{lemma1}, we see that $\Hh$ is generated over its 
center by 
$$\x^i \y^j \s^l, \quad \quad 0 \leq i,j \leq pr-1, \quad 0 \leq l 
\leq r-1.$$  
\end{proof}

\begin{corollary}
Every irreducible $\Hh$-module is finite-dimensional over $\k$.
\label{cor3}
\end{corollary}

\begin{proof}
Let $V$ be an irreducible $\Hh$-module. Then $V$ is cyclic; that is,
$V = \Hh / J$, where $J$ is a left ideal of $\Hh$. By corollary 
\ref{cor2}, it follows that $V$ is a finitely generated 
$\cent{\Hh}$-module. Hence, there exists a maximal ideal $M \subset 
\cent{\Hh}$ such that $V / MV \neq 0$. But $MV$ is a 
submodule of the irreducible $\Hh$-module $V$, so $MV = 0$. 
Thus the action of $\Hh$ in $V$ factors
through $\Hh/M \Hh$, which is a finite-dimensional algebra.
So $V=(\Hh/M \Hh)/L$ for some left ideal $L$ of $\Hh/M\Hh$. Hence $V$ is 
finite-dimensional.
\end{proof}

Thus Schur's lemma implies that 
central elements of $\Hh$ act as scalars in 
any irreducible $\Hh$-module. 

\begin{lemma} One has
$$[\hhh,\x]=\x, [\hhh,\y]=-\y, [\hhh,\s] = 0.$$
\label{lemma2}
\end{lemma}

\begin{proof}
We have
\begin{eqnarray*}
[\hhh,\x]
&=& \hhh\x - \x\hhh \\
&=& \x\y\x - \x^2 \y + \sum_{j=1}^{r-1} c_j (1 - \ert^{-j})^{-1} 
(\x \s^j - \s^j \x) \\
&=& \x [\y,\x] + \x \sum_{j=1}^{r-1} c_j \s^j \\
&=& \x [\y,\x] + \x (1 - [\y,\x]) \\
&=& \x.
\end{eqnarray*}
Similarly, 
$[\hhh,\y] = -\y.$
Finally, $\s\x\y = \ert^{-1}\x\s\y = \x\y\s$. From this it follows that 
$\s$ commutes with $\hhh$.
\end{proof}

\begin{lemma}
Let $V$ be a representation of $\Hh$. Then, the sum of the common 
eigenspaces of 
$\hhh$ and $\s$ is a subrepresentation.
\label{lemma3}
\end{lemma}

\begin{proof}
Let $v \in V$ be a common eigenvector of $\hhh$ and $\s$, with eigenvalues 
$\lambda$ and $\mu$ respectively. Then, using lemma \ref{lemma2},
\begin{eqnarray*}
\hhh (\x v) &=& \x \hhh v + \x v = (\lambda + 1) (\x v), \\ 
\hhh (\y v) &=& \y \hhh v - \y v = (\lambda - 1) (\y v), \\
\hhh (\s v) &=& \s \hhh v = \lambda (\s v), \\
\s (\x v) &=& \ert^{-1} \x\s v = \ert^{-1} \mu (\x v), \\
\s (\y v) &=& \ert \y\s v = \ert \mu (\y v), \\
\s (\s v) &=& \mu (\s v).
\end{eqnarray*}
Thus, $\x v, \y v$ and $\s v$ are also common eigenvectors of $\hhh$ and 
$\s$. The result follows.
\end{proof}

\begin{corollary}
If $V$ is an irreducible representation of $\Hh$, then $$V = \oplus_{\beta,m} 
V[\beta,m],$$ where $V[\beta,m]$ is the common eigenspace of $\hhh$ and $\s$ 
with eigenvalues $\beta$ and $\ert^m$.
\end{corollary}

\begin{corollary}
Let $\overline{\Hh}$ be the subalgebra of 
$\Hh$ generated by $\cent{\Hh}, \hhh$ and 
$\s$. Then $\overline{\Hh}$ is commutative.
\label{cornew}
\end{corollary}

\begin{proof}
Follows from lemma \ref{lemma2}.
\end{proof}

\begin{lemma} One has
$\x^m \y^n \s^l \in \overline{\Hh}$ for all $m,n,l$ such that $pr$
divides $m-n$.
\label{lemmanew}
\end{lemma}

\begin{proof}
For all $m \geq 1$, we have
\begin{eqnarray*}
\x^m \y^m &=& \x^{m-1} (\x \y) \y^{m-1} \\
&=& \x^{m-1} (\hhh + f_0(\s)) \y^{m-1} \quad \mbox{for
some } f_0(\s) \in
\k[\s] \\
&=& f_1(\s) \x^{m-1} \y^{m-1} + \x^{m-1} \hhh \y^{m-1}
\quad \mbox{for
some } f_1(\s) \in \k[\s] \\
&=& f_1(\s) \x^{m-1} \y^{m-1} + \x^{m-1} \y^{m-1} \hhh -
(m-1) \x^{m-1} \y^{m-1}.   
\end{eqnarray*}
Using induction, it follows that $\x^m \y^m \in \overline{\Hh}$
for all $m \geq 0$.
Using lemma \ref{lemma1}, we see that $\x^m \y^n \s^l \in \overline{\Hh}$ whenever
$pr$ divides $m-n$. 
\end{proof}

\begin{corollary}
Let $V$ be an irreducible $\Hh$-module. Then $\dim V[\beta,m] \leq 1$ for each 
$\beta \in \k, m \in \Z / r\Z$. Moreover, if $V[\beta,m] \neq 0$ and 
$V[\beta', m'] \neq 0$, then $\beta-\beta'$ is an element of 
the prime field $\pf \subset \k$, and thus $\dim V \leq pr$.
\label{corlast}
\end{corollary}

\begin{proof}
Let $v \in V$ be nonzero. Since $V$ is irreducible, we have $\Hh v = V$.
By lemma \ref{PBW}, we know that $V$ is generated over $\k$ by
\begin{equation*}
\x^i \y^j \s^l v, \quad 0 \leq i,j, 0 \leq l \leq r-1.
\end{equation*}
By lemma \ref{lemma1}, it follows that $V$ is generated 
over $\k$ by  
\begin{equation*}
\x^i \y^j \s^l v, \quad 0 \leq i,j \leq pr-1, 0 \leq l \leq r-1.
\end{equation*}

Now suppose $V[\beta,m] \neq 0$, and let
$v \in V[\beta,m]$ be nonzero. Then $\x^i \y^j \s^l v \in V[\beta+i-j, 
m-i+j]$ for all $i,j,l$. This proves the second part of the corollary, and 
also implies that 
$V[\beta,m]$ is generated by $\x^i \y^i \s^l v$ with $0 \leq i \leq 
pr-1$. Thus $V[\beta, m] \subset \overline{\Hh} v$, by lemma \ref{lemmanew}.
But it is easy to see that $\overline{\Hh} v \subset V[\beta,m]$.
Thus $\overline{\Hh} v = V[\beta, m]$. Since this is true for all nonzero $v \in 
V[\beta,m]$, we conclude that $V[\beta,m]$ is an irreducible 
representation of $\overline{\Hh}$. 
Since $\overline{\Hh}$ is commutative, we see that 
$V[\beta,m]$ is one dimensional. This proves the first part 
of the corollary.
\end{proof}

\begin{proof}[Proof of proposition \ref{prop1}]
Let $V$ be an irreducible representation of $\Hh$.

Let $v \in V[\beta,m]$ be nonzero. Then
$\x^i v \in V[\beta+i, m-i]$ for all $i$. By lemma 
\ref{lemma1}, we know that $\x^{pr}, \y^{pr}$ act on $V$ by multiplication
by scalars $a,b$ respectively. Assume $a \neq 0$. Then $\x^{pr} v 
= av \neq 0$, so $\x^i v \neq 0$ for all $0 \leq i \leq pr-1$, and
thus $V[\beta+i,m-i] \neq 0$ for $0 \leq i \leq pr-1$. Applying
corollary \ref{corlast}, and using the fact that $p,r$ are 
relatively prime, we conclude that $\dim V = pr$ and 
$\dim V[\beta',m'] = 1$ for $\beta' = \beta, \beta+1, \ldots, 
\beta+(p-1), 0 \leq m' \leq r-1$. We may now assume that $m = 0$.  

Now let $v_i \defeq \x^i v$ for $i = 0, 1, \ldots, pr$. In 
particular, $v_0 = v$ and $v_{pr} = a v$. Then the $v_i, 0 \leq i 
\leq pr-1$, form a basis of $V$. We have $v_i \in V[\beta+i, -i]$, so
$\s v_i = \ert^{-i} v_i$. For each $i$, $\y v_i \in V[\beta+i-1, -i+1]$, 
so we let $\y v_i = \mu_i v_{i-1}$ for $1 \leq i \leq pr$. Then, for $1 
\leq i \leq pr$, we have 
\begin{eqnarray*}
(\beta+i) v_i &=& \hhh v_i \\
&=& \left(\x\y - \sum_{j=1}^{r-1} c_j 
	(1 - \ert^{-j})^{-1} \s^j \right) v_i \\
&=& \left(\mu_i - \sum_{j=1}^{r-1} c_j (1 - \ert^{-j})^{-1} 
	\ert^{-ij} \right) v_i,
\end{eqnarray*}
so $$\mu_i = \beta + i + \sum_{j=1}^{r-1} c_j (1 - \ert^{-j})^{-1}
        \ert^{-ij},$$ which proves uniqueness.
Now it is easy to check that $V_{\beta,a,b} = \span \{v_i, 0 \leq i \leq 
pr-1\}$,
\begin{eqnarray*}
\x v_i &=& v_{i+1}, \quad 0 \leq i \leq pr-1 \\
\y v_i &=& \mu_i v_{i-1}, \quad 1 \leq i \leq pr \\
\s v_i &=& \ert^{-i} v_i, \quad 0 \leq i \leq pr \\
v_{pr} &=& a v_0
\end{eqnarray*}
defines a representation of $\Hh$. Furthermore, the representation is 
irreducible, since any nonzero subrepresentation of $V_{\beta,a,b}$ would 
have dimension $pr$, by the arguments contained in this proof.
We note that $$ab v_0 = b v_{pr} = \y^{pr} v_{pr} = 
\left(\prod_{i=1}^{pr} \mu_i \right) v_0.$$
Therefore, 
\begin{eqnarray*}
ab &=& \prod_{i=1}^{pr} \left( \beta + i + \sum_{j=1}^{r-1} c_j 
	(1 - \ert^{-j})^{-1} \ert^{-ij} \right) \\
&=& \prod_{m=0}^{r-1} \prod_{i=0}^{p-1} \left(\beta + i + 
	\sum_{j=1}^{r-1} c_j (1 - \ert^{-j})^{-1} \ert^{mj} \right) \\
&=& \prod_{m=0}^{r-1} f \left(\beta + \sum_{j=1}^{r-1} c_j (1 - 
	\ert^{-j})^{-1} \ert^{mj} \right),
\end{eqnarray*}
where $f(z) = z^p - z$.
It remains to show that $V_{\beta,a,b} \isom V_{\beta',a,b}$ if and only 
if $\beta - \beta' \in \F_p$. But if $\beta - \beta' \in \F_p$, then the 
common eigenvalues of $\hhh$ and $\s$ for $V_{\beta,a,b}$ and 
$V_{\beta',a,b}$ are
the same, so $V_{\beta,a,b} \isom V_{\beta',a,b}$ by the arguments in this 
proof. If $\beta - \beta' \notin \F_p$, then the
common eigenvalues of $\hhh$ and $\s$ for $V_{\beta,a,b}$ and 
$V_{\beta',a,b}$ 
are different, so $V_{\beta,a,b}$ and $V_{\beta'a,b}$ cannot be 
isomorphic.

A similar argument works in the case $a = 0, b \neq 0$.
\end{proof}

\begin{proof}[Proof of proposition \ref{prop2}]
Let $V$ be an irreducible representation of $\Hh$ with $\x^{pr} = \y^{pr} 
= 0$.
Let $v$ be a common eigenvector of $\hhh$ and $\s$. Since $\x^{pr} V = 0$ 
and $\y^{pr} V = 
0$, we can find minimal $d, D$ such that $\x^D v = \y^d v = 0$. Let
$$v_{-i} = \y^i v, 1 \leq i < d; \quad v_0 = v; \quad v_i = \x^i v, 1 \leq 
i < D.$$ Since $V$ is irreducible, we know that the $v_i$ form a basis of 
$V$ and that $\x v_{-i} \neq 0$ for $1 \leq i < d$, $\y v_i \neq 0$ for $1 
\leq i < D$. Thus, by scaling and relabeling, and using corollary 
\ref{corlast}, we may assume that $V$ has a basis $v_0, v_1, \ldots, 
v_{D-1}$ with 
\begin{eqnarray*}
\x v_i &=& v_{i+1}, \quad \mbox{for } 0 \leq i \leq D-2; \\
\x v_{D-1} &=& 0; \\
\y v_i &=& \mu_i v_{i-1} \quad \mbox{for } 1 \leq i \leq D-1; \\
\y v_0 &=& 0; \\
\s v_i &=& \ert^{m-i} v_i \quad \mbox{for } 0 \leq i \leq D-1.
\end{eqnarray*}
In order for this to be a representation, we need
$$(\y\x - \x\y) v_i = = \left(1 - \sum_{j=1}^{r-1} c_j \s^j \right) v_i.$$
Writing $\mu_D = \mu_0 = 0$, this gives
\begin{equation}
\mu_{i+1} - \mu_i = 1 - \sum_{j=1}^{r-1} c_j \ert^{(m-i)j} \quad 
\mbox{for } 0 \leq i \leq D-1.
\label{mu}
\end{equation}
Adding up \eqref{mu} for $0 \leq i \leq k-1$ gives
$$\mu_k = k - \sum_{i=0}^{k-1} \sum_{j=1}^{r-1} c_j \ert^{(m-i)j}.$$
In particular,
$$0 = \mu_D = D - \sum_{i=0}^{D-1} \sum_{j=1}^{r-1} c_j \ert^{(m-i)j}.$$
As we saw earlier, in order for the representation to be irreducible,
we need $\mu_k \neq 0$ for $1 \leq k \leq D-1$, which translates as
$$0 \neq k - \sum_{i=0}^{k-1} \sum_{j=1}^{r-1} c_j \ert^{(m-i)j}, 
\quad 1 \leq k \leq D-1.$$
\end{proof}

\begin{corollary}[PBW for $\Hh$]
The elements $$\x^i \y^j \s^l, \quad \quad 0 \leq i,j, \quad 0 \leq l
\leq r-1$$ form a basis of $\Hh$ over $\k$.
\label{hPBW}
\end{corollary}

\begin{proof}
This follows because the elements 
$\x^i \y^j \s^l$ act as linearly independent operators on 
$\oplus_{\beta,a,b}V_{\beta,a,b}$. 
\end{proof}

\section{Proof of propositions \ref{0prop1} and \ref{0prop2}}

In this section, several proofs will be omitted because they are essentially
identical to the corresponding proofs in the previous section.

\begin{lemma}[PBW for $\Hh'$, easy direction]
The elements $$\x^i \y^j \s^l, \quad \quad 0 \leq i,j, \quad 0 \leq l
\leq r-1$$ span $\Hh'$ over $\k$.
\label{0PBW}
\end{lemma}   

\begin{lemma}
$\x^{r}$ and $\y^{r}$ belong to the center $\cent{\Hh'}$ of $\Hh'$.
\label{0lemma1}
\end{lemma}

\begin{proof} 
First, let us show that $\x^{r} \in \cent{\Hh'}$.
It is enough to show that $[\x^{r},\y] = [\x^{r},\s] = 0$.
But $\x^r \s = \ert^r \s \x^r = \s \x^r$, so $[\x^r, \s] = 0$.
Now let $g(\s) = - \sum_{j=1}^{r-1} c_j \s^j.$ Then,
\begin{eqnarray*}
\y\x^r &=& (\y\x)\x^{r-1} \\
&=& (\x\y + g(\s)) \x^{r-1} \\
&=& \x (\y \x) \x^{r-2} + g(\s) \x^{r-1} \\
&=& \cdots \\
&=& \x^r\y + g(\s) \x^{r-1} + \x g(\s) \x^{r-2} + \cdots + \x^{r-1} g(\s) \\
&=& \x^r\y + \sum_{j=1}^{r-1} c_j \left( \s^j \x^{r-1} + \x
\s^j \x^{r-2} + \cdots + \x^{r-1} \s^j \right) \\
&=& \x^r\y + \sum_{j=1}^{r-1} c_j (1 + \ert^j + \cdots +
\ert^{j(r-1)}) \s^j \x^{r-1} \\
&=& \x^r\y,
\end{eqnarray*}
since $1 + z + \cdots + z^{r-1} = 0$ whenever $z \neq 1$ is a
$r$-th root of unity. Thus $[\y, \x^r] = 0$.

This completes the proof that $\x^{r} \in \cent{\Hh'}.$ The proof that
$\y^{r} \in \cent{\Hh'}$ is similar.
\end{proof}

\begin{corollary}
$\Hh'$ is finitely generated as a module over its center.
\label{0cor2}
\end{corollary}

\begin{corollary}
Every irreducible $\Hh'$-module is finite-dimensional over $\k$.
\label{0cor3}
\end{corollary}

In particular, central elements of $\Hh'$ act as scalars in its irrreducible 
representations.

\begin{lemma}
One has $\hhh \in \cent{\Hh'}.$
\label{0lemma2}
\end{lemma}

\begin{proof}
We have
\begin{eqnarray*}  
[\hhh,\x]  
&=& \hhh\x - \x\hhh \\
&=& \x\y\x - \x^2 \y + \sum_{j=1}^{r-1} c_j (1 - \ert^{-j})^{-1}
(\x \s^j - \s^j \x) \\
&=& \x [\y,\x] + \x \sum_{j=1}^{r-1} c_j \s^j \\
&=& \x [\y,\x] + \x (- [\y,\x]) \\
&=& 0.
\end{eqnarray*}
Similarly,   
$[\hhh,\y]= 0.$
Finally, $\s\x\y = \ert^{-1}\x\s\y = \x\y\s$, so $[\hhh, \s] = 0$.
\end{proof}

Let $\overline{\Hh'}$ be the subalgebra of 
$\Hh'$ generated by $\cent{\Hh'}$ and
$\s$. Obviously, $\overline{\Hh'}$ is commutative.

\begin{lemma}
$\x^m \y^n \s^l \in \overline{\Hh'}$ for all $m,n,l$ such that $r$
divides $m-n$.
\label{0lemmanew}
\end{lemma}

\begin{proof}  
For all $m \geq 1$, we have
\begin{eqnarray*}
\x^m \y^m &=& \x^{m-1} (\x \y) \y^{m-1} \\
&=& \x^{m-1} (\hhh + f_0(\s)) \y^{m-1} \quad \mbox{for
some } f_0(\s) \in
\k[\s] \\
&=& f_1(\s) \x^{m-1} \y^{m-1} + \x^{m-1} \y^{m-1} \hhh
\quad \mbox{for
some } f_1(\s) \in \k[\s].
\end{eqnarray*}
Using induction, it follows that $\x^m \y^m \in \overline{\Hh'}$
for all $m \geq 0$.
Using lemma \ref{0lemma1}, we see that $\x^m \y^n \s^l \in \overline{\Hh'}$ whenever
$r$ divides $m-n$.   
\end{proof}

If $V$ is an irreducible representation of $\Hh'$, then $$V = \oplus_{m}
V[m],$$ where $V[m]$ is the eigenspace of $\s$
with eigenvalue $\ert^m$.

\begin{corollary}
Let $V$ be an irreducible $\Hh'$-module. Then $\dim V[m]
\leq 1$ for each
$m \in \Z / r\Z$.
Thus $\dim V \leq r$.
\label{0corlast}
\end{corollary}

\begin{proof}
Let $v \in V$ be nonzero. Since $V$ is irreducible, we have $\Hh' v = V$.
By lemma \ref{0PBW}, we know that $V$ is generated over $\k$ by
\begin{equation*} 
\x^i \y^j \s^l v, \quad 0 \leq i,j, 0 \leq l \leq r-1.
\end{equation*}
By lemma \ref{0lemma1}, it follows that $V$ is generated
over $\k$ by
\begin{equation*}
\x^i \y^j \s^l v, \quad 0 \leq i,j \leq r-1, 0 \leq l \leq r-1.
\end{equation*}

Now suppose $V[m] \neq 0$, and let
$v \in V[m]$ be nonzero. Then $\x^i \y^j \s^l v \in V[m-i+j]$ for all
$i,j,l$. This implies that
$V[m]$ is generated over $\k$ by $\x^i \y^i \s^l v$ with $0 \leq i \leq
r-1$. Thus $V[m] \subset \overline{\Hh'} v$, by lemma \ref{0lemmanew}.
But it is easy to see that $\overline{\Hh'} v \subset V[m]$.
Thus $\overline{\Hh'} v = V[m]$. Since this is true for all nonzero $v \in
V[m]$, we conclude that $V[m]$ is an irreducible
representation of $\overline{\Hh'}$, hence one dimensional.
This proves the
corollary.
\end{proof}

\begin{proof}[Proof of proposition \ref{0prop1}]
Let $V$ be an irreducible representation of $\Hh'$.

Let $v \in V[m]$ be nonzero. Then
$\x^i v \in V[m-i]$ for all $i$. By lemmas \ref{0lemma1}
and \ref{0lemma2}, we know that $\x^{r}, \y^r$ and $\hhh$ act on $V$ by
multiplication
by scalars, say $a,b,\beta$ respectively. Assume $a \neq 0$. Then $\x^{r} 
v = av \neq 0$, so $\x^i v \neq 0$ for all $0 \leq i \leq r-1$, and
thus $V[m-i] \neq 0$ for $0 \leq i \leq r-1$. Applying
corollary \ref{0corlast}, we conclude that $\dim V = r$ and $\dim V[m'] = 
1$ for $0 \leq m' \leq r-1$. We may now assume that $m = 0$.

Now let $v_i \defeq \x^i v$ for $i = 0, 1, \ldots, r$. In
particular, $v_0 = v$ and $v_{r} = a v$. Then the $v_i, 0 \leq i
\leq r-1$, form a basis of $V$. We have $v_i \in V[-i]$, so
$\s v_i = \ert^{-i} v_i$. For each $i$, $\y v_i \in V[-i+1]$,
so we let $\y v_i = \mu_i v_{i-1}$ for $1 \leq i \leq r$. Then, for $1
\leq i \leq r$, we have
\begin{eqnarray*}
\beta v_i &=& \hhh v_i \\
&=& \left(\x\y - \sum_{j=1}^{r-1} c_j
        (1 - \ert^{-j})^{-1} \s^j \right) v_i \\
&=& \left(\mu_i - \sum_{j=1}^{r-1} c_j (1 - \ert^{-j})^{-1}
        \ert^{-ij} \right) v_i,
\end{eqnarray*}
so $$\mu_i = \beta + \sum_{j=1}^{r-1} c_j (1 - \ert^{-j})^{-1}
        \ert^{-ij},$$ which proves uniqueness.
Now it is easy to check that $V_{\beta,a,b} = \span \{v_i, 0 \leq i \leq
r-1\}$,
\begin{eqnarray*}
\x v_i &=& v_{i+1}, \quad 0 \leq i \leq r-1 \\
\y v_i &=& \mu_i v_{i-1}, \quad 1 \leq i \leq r \\
\s v_i &=& \ert^{-i} v_i, \quad 0 \leq i \leq r \\
v_{r} &=& a v_0
\end{eqnarray*}
defines a representation of $\Hh'$. Furthermore, the representation is
irreducible, since any nonzero subrepresentation of $V_{\beta,a,b}$ would
have dimension $r$, by the arguments contained in this proof.
We note that $$ab v_0 = b v_{r} = \y^{r} v_{r} =
\left(\prod_{i=1}^{r} \mu_i \right) v_0.$$
Therefore,
\begin{equation*}
ab = \prod_{i=1}^{r} \left( \beta + \sum_{j=1}^{r-1} c_j
        (1 - \ert^{-j})^{-1} \ert^{-ij} \right).
\end{equation*}

A similar argument works in the case $b \neq 0$.
\end{proof}

\begin{proof}[Proof of proposition \ref{0prop2}]
Let $V$ be an irreducible representation of $\Hh'$ with $\x^{r} = \y^{r}
= 0$.
Let $v$ be an eigenvector of $\s$. Since $\x^{r} V = 0$
and $\y^{r} V = 0$, we can find minimal $d, D$ such that $\x^D v = \y^d v = 0$. 
Let $$v_{-i} = \y^i v, 1 \leq i < d; \quad v_0 = v; \quad v_i = \x^i v, 1 \leq
i < D.$$ Since $V$ is irreducible, we know that the $v_i$ form a basis of
$V$ and that $\x v_{-i} \neq 0$ for $1 \leq i < d$, $\y v_i \neq 0$ for $1
\leq i < D$. Thus, by scaling and relabeling, and using corollary
\ref{0corlast}, we may assume that $V$ has a basis $v_0, v_1, \ldots,
v_{D-1}$ with
\begin{eqnarray*}
\x v_i &=& v_{i+1}, \quad \mbox{for } 0 \leq i \leq D-2; \\  
\x v_{D-1} &=& 0; \\
\y v_i &=& \mu_i v_{i-1} \quad \mbox{for } 1 \leq i \leq D-1; \\
\y v_0 &=& 0; \\
\s v_i &=& \ert^{m-i} v_i \quad \mbox{for } 0 \leq i \leq D-1.
\end{eqnarray*}
In order for this to be a representation, we need
$$(\y\x - \x\y) v_i = - \sum_{j=1}^{r-1} c_j \s^j v_i.$$
Writing $\mu_D = \mu_0 = 0$, this gives
\begin{equation}
\mu_{i+1} - \mu_i = - \sum_{j=1}^{r-1} c_j \ert^{(m-i)j} \quad
\mbox{for } 0 \leq i \leq D-1.
\label{0mu}
\end{equation}
Adding up \eqref{0mu} for $0 \leq i \leq t-1$ gives
$$\mu_k = - \sum_{i=0}^{k-1} \sum_{j=1}^{r-1} c_j \ert^{(m-i)j}.$$
In particular,
$$0 = \mu_D = - \sum_{i=0}^{D-1} \sum_{j=1}^{r-1} c_j \ert^{(m-i)j}.$$
As we saw earlier, in order for the representation to be irreducible,
we need $\mu_k \neq 0$ for $1 \leq k \leq D-1$, which translates as
$$0 \neq \sum_{i=0}^{k-1} \sum_{j=1}^{r-1} c_j \ert^{(m-i)j}, 
\quad 1 \leq k \leq D-1.$$
\end{proof}

\begin{corollary}[PBW for $\Hh'$]                
The elements $$\x^i \y^j \s^l, \quad \quad 0 \leq i,j, 
\quad 0 \leq l \leq r-1$$ form a basis of $\Hh'$ over $\k$.
\label{0hPBW}
\end{corollary}

\begin{proof}
This follows because the 
elements $\x^i \y^j \s^l$ act as linearly independent operators 
on $\oplus_{\beta,a,b}V_{\beta,a,b}$.
\end{proof}

\bibliography{prob6}
\bibliographystyle{alpha}

\end{document}